# ALMOST GLOBAL EXISTENCE FOR SOME SEMILINEAR WAVE EQUATIONS

MARKUS KEEL, HART F. SMITH, AND CHRISTOPHER D. SOGGE

*Tom Wolff, in memoriam*

1. **Introduction.**

This article establishes the almost global existence of solutions of three-dimensional quadratically semilinear wave equations with the use of only the classical invariance of the equations under translations and spatial rotation. Using these techniques we can handle semilinear wave equations in Minkowski space or semilinear Dirichlet-wave equations in the exterior of a nontrapping obstacle.

Our results and approach are related to previous work in the non-obstacle case. In particular, in [2], almost global existence was shown for general, quadratic *quasilinear* wave equations in Minkowski space, using fully the invariance properties of the wave operator under translations, spatial rotations, scaling, and Lorentz boosts. Subsequently, in [5], the almost global existence of certain divergence-form quadratic quasilinear wave equations was established in Minkowski space, using only the invariance of the wave operator under translations, spatial rotations, and changes of scale. Further techniques and results concerning equations from the theory of elasticity were presented in [9], [10].

The main difference between our approach and earlier work is that it only exploits the $O(|x|^{-1})$ decay of solutions of wave equations with compactly supported data, as opposed to the stronger $O(t^{-1})$ decay. Here of course $x = (x_1, x_2, x_3)$ is the spatial component and $t$ the time component of a space-time vector $(t, x) \in \mathbb{R}_+ \times \mathbb{R}^3$. Establishing $O(|x|^{-1})$ decay is considerably easier and can be achieved using only the invariance with respect to translations and spatial rotations. Fortunately, as we shall see, one can set up an iteration argument that only requires this weaker decay to obtain almost global existence. Another advantage to this approach for obstacle problems is that the translation and spatial rotational vector fields essentially preserve the Dirichlet conditions. The homogeneous vector fields used in [2] do not have this property for any obstacle, and consequently it seems difficult to use them for nonlinear obstacle problems.

Let us now describe more precisely the equations that we shall consider. If $\mathcal{K} \subset \mathbb{R}^3$ is a smooth compact nontrapping obstacle, we shall consider semilinear systems of the

The authors were supported in part by the NSF.

The third author would like to thank J.-M. Delort for helpful conversations regarding this paper. He would also like to thank Paris Nord for their hospitality during a recent visit where part of this research was carried out.

form

(1.1)
$$\begin{cases} \Box u = Q(u'), \quad (t,x) \in \mathbb{R}_+ \times \mathbb{R}^3 \backslash \mathcal{K} \\ u(t,\cdot)|_\mathcal{K} = 0 \\ u(0,\cdot) = f, \ \partial_t u(0,\cdot) = g. \end{cases}$$

Here
$$\Box = \partial_t^2 - \Delta$$
is the D'Alembertian, with $\Delta = \partial_1^2 + \partial_2^2 + \partial_3^2$ being the standard Laplacian. Also, $Q$ is a constant coefficient quadratic form in $u' = (\partial_t u, \nabla_x u)$.

In the non-obstacle case we shall obtain almost global existence for equations of the form

(1.2)
$$\begin{cases} \Box u = Q(u'), (t,x) \in \mathbb{R}_+ \times \mathbb{R}^3 \\ u(0,\cdot) = f, \ \partial_t u(0,\cdot) = g. \end{cases}$$

In order to solve (1.1) we must also assume that the data satisfies the relevant compatibility conditions. Since these are well known (see e.g., [3]), we shall describe them briefly. To do so we first let $J_k u = \{\partial_x^\alpha u : 0 \le |\alpha| \le k\}$ denote the collection of all spatial derivatives of $u$ of order up to $k$. Then if $m$ is fixed and if $u$ is a formal $H^m$ solution of (1.1) we can write $\partial_t^k u(0,\cdot) = \psi_k(J_k f, J_{k-1} g)$, $0 \le k \le m$, for certain compatibility functions $\psi_k$ which depend on the nonlinear term $Q$ as well as $J_k f$ and $J_{k-1} g$. Having done this, the compatibility condition for (1.1) with $(f,g) \in H^m \times H^{m-1}$ is just the requirement that the $\psi_k$, vanish on $\partial \mathcal{K}$ when $0 \le k \le m-1$. Additionally, we shall say that $(f,g) \in C^\infty$ satisfy the compatibility conditions to infinite order if this condition holds for all $m$.

If $\{\Omega\}$ denotes the collection of vector fields $x_i \partial_j - x_j \partial_i$, $1 \le i < j \le 3$, then we can now state our main result.

**Theorem 1.1.** *Let $\mathcal{K}$ be a smooth compact nontrapping obstacle and assume that $Q(u')$ is above. Assume further that $(f,g) \in C^\infty(\mathbb{R}^3 \backslash \mathcal{K})$ satisfies the compatibility conditions to infinite order. Then there are constants $c, \varepsilon_0 > 0$ so that if $\varepsilon \le \varepsilon_0$ and*

(1.3)
$$\sum_{|\alpha|+j\le 10} \|\partial_x^j \Omega^\alpha f\|_{L^2(\mathbb{R}^3\backslash\mathcal{K})} + \sum_{|\alpha|+j\le 9} \|\partial_x^j \Omega^\alpha g\|_{L^2(\mathbb{R}^3\backslash\mathcal{K})} \le \varepsilon,$$

*then (1.1) has a unique solution $u \in C^\infty([0,T_\varepsilon] \times \mathbb{R}^3 \backslash \mathcal{K})$, with*

(1.4)
$$T_\varepsilon = \exp(c/\varepsilon).$$

We shall actually establish existence of limited regularity almost global solutions $u$ for data $(f,g) \in H^9 \times H^8$ satisfying the relevant compatibility conditions and smallness assumptions (1.3). The fact then that $u$ must be smooth if $f$ and $g$ are smooth and satisfy the compatibility conditions of infinite order follows from standard local existence theorems (see §9, [3]).

We should point out that results similar to those in Theorem 1.1 were announced in Datti [1], but there appears to be a gap in the argument which has not been repaired.

As we pointed out before, we can also give a new proof of the semilinear version of the almost global existence theorem of John and Klainerman [2]:

**Theorem 1.2.** *Assume that $Q(u')$ is above, and that $(f,g) \in C^\infty(\mathbb{R}^3)$. Then there is a constant $c > 0$ so that if $\varepsilon > 0$ is small and*

$$\sum_{|\alpha|+j \leq 10} \|\partial_x^j \Omega^\alpha f\|_{L^2(\mathbb{R}^3)} + \sum_{|\alpha|+j \leq 9} \|\partial_x^j \Omega^\alpha g\|_{L^2(\mathbb{R}^3)} \leq \varepsilon. \tag{1.5}$$

*then (1.2) has a unique solution $u \in C^\infty([0, T_\varepsilon] \times \mathbb{R}^3)$, with*

$$T_\varepsilon = \exp(c/\varepsilon). \tag{1.6}$$

This paper is organized as follows. In the next section we shall collect the main estimates that we require. In §3 we shall give the simple proof of Theorem 1.2. After that we shall discuss the straightforward modifications that are needed to obtain Theorem 1.1.

In a later paper we hope to use similar techniques to handle quasilinear equations.

2. **Main estimates.**

In what follows $r = |x|$. A key estimate that we require is the following radial decay estimate.

**Proposition 2.1.** *Suppose that $v$ solves the wave equation $\Box v = G$ in $\mathbb{R}_+ \times \mathbb{R}^3$. Then there is a uniform constant $C$ so that*

$$(\ln(2+t))^{-1/2}\|(1+r)^{-1/2}v'\|_{L^2(\{(s,x):\, 0\leq s \leq t\})} \tag{2.1}$$
$$\leq C\|v'(0,\,\cdot\,)\|_{L^2(\mathbb{R}^3)} + C\int_0^t \|G(s,\cdot)\|_{L^2(\mathbb{R}^3)}\,ds.$$

This will follow as a consequence of the following

**Lemma 2.2.** *Let $v$ be as above. Then*

$$\|v'\|_{L^2(\{(s,x):\, 0\leq s\leq t,\, |x|<1\})} \leq C\|v'(0,\,\cdot\,)\|_{L^2(\mathbb{R}^3)} + C\int_0^t \|G(s,\cdot)\|_{L^2(\mathbb{R}^3)}\,ds. \tag{2.2}$$

*Also,*

$$(1+t)^{-1/2}\|v'\|_{L^2(\{(s,x):\, 0\leq s\leq t\})} \leq C\|v'(0,\,\cdot\,)\|_{L^2(\mathbb{R}^3)} + C\int_0^t \|G(s,\cdot)\|_{L^2(\mathbb{R}^3)}\,ds. \tag{2.3}$$

To see that Lemma 2.2 implies Proposition 2.1, note first that on the set $r > t$, (2.1) follows directly from (2.3). On the set $r \leq t$, make a dyadic decompositon in $r$ and for each piece use the following estimate, which follows from (2.2) and a scaling argument,

$$\|r^{-1/2}v'\|_{L^2(\{(s,x):\, 0\leq s\leq t,\, |x|\in [R,2R]\})} \leq C\|v'(0,\,\cdot\,)\|_{L^2(\mathbb{R}^3)} + C\int_0^t \|G(s,\cdot)\|_{L^2(\mathbb{R}^3)}\,ds. \tag{2.4}$$

By squaring the left side of (2.1) and then summing over the dyadic pieces with $1 < r < t$, (2.4) together with (2.2) give (2.1).

**Proof of Lemma 2.2:** Note that (2.3) is a direct consequence of the energy inequality and Young's inequality in the $t$ variable.

To prove the estimate (2.2), we will assume that $v$ has vanishing initial data. A straightforward modification of the argument will handle the case where the Cauchy data is non-zero.

We shall use sharp Huygen's principle. Split $G(t,x) = \sum_{j\geq 1} G_j(s,x)$ where $G_j(s,x) = G(s,x)$ if $s + |x| \in (j-1, j]$ and $0$ otherwise. Then if $(s,x) \in [j, j+1] \times \{x : |x| < 1\}$ it follows that $v = v_j$ where $v_j$ solves the wave equation $\Box v_j = \sum_{|j-k|\leq 10} G_k$, with vanishing initial data because of our assumptions. On account of this the energy inequality gives

$$\|v'\|_{L^2(\{(s,x)\colon s\in[j,j+1], |x|<1\})} \leq \sum_{|j-k|\leq 10} \int_0^t \|G_k(s,\cdot)\|_{L^2(\mathbb{R}^3)}\, ds,$$

assuming that $j+1 < t$. The bound now follows from the fact that

$$\sum_k \left(\int_0^t \|G_k(s,\cdot)\|_{L^2(\mathbb{R}^3)}\, ds\right)^2 \leq C \left(\int_0^t \|G(s,\cdot)\|_{L^2(\mathbb{R}^3)} ds\right)^2,$$

because of the disjoint supports of the functions. $\square$

For later use, note that this proof also implies that

(2.5) $$\|v\|_{L^2(\{(s,x)\colon 0\leq s\leq t, |x|<1\})} \leq C\|v'(0,\cdot)\|_{L^2(\mathbb{R}^3)} + C\int_0^t \|G(s,\cdot)\|_{L^2(\mathbb{R}^3)}\, ds,$$

since

$$\|v(s,\cdot)\|_{L^2(|x|<1)} \leq \|v(s,\cdot)\|_{L^6(|x|<1)} \leq \int_0^s \|G(\tau,\cdot)\|_{L^2(\mathbb{R}^3)}\, d\tau.$$

Let $\{Z\} = \{\partial_x, \partial_t, \Omega\}$, where as before the $\{\Omega\}$ are the $\mathbb{R}^3$ rotational vector fields. Then since all of the $Z$ commute with $\Box$, Proposition 2.1 and the energy inequality imply the following

**Theorem 2.3.** *Let $v$ be as above. Then for any $N = 0, 1, 2, \ldots$*

(2.6) $$\sum_{|\alpha|\leq N}\left(\|Z^\alpha v'(t,\cdot)\|_2 + (\ln(2+t))^{-1/2}\|(1+r)^{-1/2}Z^\alpha v'\|_{L^2(\{(s,x)\colon 0\leq s\leq t\})}\right)$$

$$\leq C\sum_{|\alpha|\leq N}\|(Z^\alpha v)'(0,\cdot)\|_2 + C\sum_{|\alpha|\leq N}\int_0^t \|Z^\alpha G(s,\cdot)\|_2 ds.$$

We shall also use the following now standard weighted Sobolev estimate. (See [4].)

**Lemma 2.4.** *Suppose that $h \in C^\infty(\mathbb{R}^3)$. Then for $R > 1$*

(2.7) $$\|h\|_{L^\infty(R/2\leq |x|\leq R)} \leq CR^{-1}\sum_{|\alpha|\leq 2}\|Z^\alpha h\|_{L^2(R/4\leq |x|\leq 2R)}.$$

**Proof:** Suppose that $|x| = r_0 \in (R/2, R)$. Then by Sobolev's lemma for $\mathbb{R} \times S^2$ we have

$$|h(x)| \leq C\sum_{|\alpha|+j\leq 2}\left(\int_{r_0-1/4}^{r_0+1/4}\int_{S^3}|\partial_r^j \Omega^\alpha h(r\omega)|^2 dr d\omega\right)^{1/2}.$$

On the other hand since the volume element in $\mathbb{R}^3$ is a constant times $r^2 dr d\omega$, the right hand side is dominated by

$$R^{-1} \sum_{|\alpha| \leq 2} \|Z^\alpha h\|_{L^2(R/4 \leq |x| \leq 2R)},$$

which finishes the proof. $\square$

## 3. Almost global existence in Minkowski space.

In this section we shall prove Theorem 1.2 using Theorem 2.3 and Lemma 2.4.

We shall solve (1.2) by an iteration argument and show that the solution satisfies

$$(3.1) \quad \sum_{|\alpha| \leq 10} \left( \|Z^\alpha u'(t, \cdot)\|_{L^2(\mathbb{R}^3)} + (\ln(2+t))^{-1/2} \|(1+r)^{-1/2} Z^\alpha u'\|_{L^2(\{(s,x):\, 0 \leq s \leq t\})} \right) \leq C\varepsilon,$$

for $0 \leq t \leq T_\varepsilon$ with $C$ being a uniform constant.

**Proof of Theorem 1.2:** If $u_{-1} = 0$, we then define $u_k$, $k = 0, 1, 2, \ldots$ inductively by letting $u_k$ solve

$$(3.2) \quad \begin{cases} \Box u_k(t,x) = Q(u'_{k-1}(t,x)), & (t,x) \in [0, T_*] \times \mathbb{R}^3 \\ u(0, \cdot) = f, \quad \partial_t u(0, \cdot) = g. \end{cases}$$

Let

$$M_k(T) = \sup_{0 \leq t \leq T} \sum_{|\alpha| \leq 10} \Big( \|Z^\alpha u'_k(t, \cdot)\|_{L^2(\mathbb{R}^3)}$$
$$+ (\ln(2+t))^{-1/2} \|(1+r)^{-1/2} Z^\alpha u'_k\|_{L^2(\{(s,x):\, 0 \leq s \leq T\})} \Big)$$

Clearly, then (1.5) and (2.6) imply that there is a uniform constant $C_0$ so that

$$M_0(T) \leq C_0 \varepsilon,$$

for any $T$, where $C_0$ is the constant in (2.6). We claim that if $\varepsilon < \varepsilon_0$ is sufficiently small and if the $c$ occurring in the definition of $T_\varepsilon$ is sufficiently small then for every $k = 1, 2, 3, \ldots$

$$(3.3) \quad M_k(T_\varepsilon) \leq 2C_0 \varepsilon.$$

Let us prove this inductively. So we shall assume that the bound holds for $k - 1$ and then prove it for $k$. By (2.6)

$$(3.4) \quad M_k(T_\varepsilon) \leq C_0 \varepsilon + C \sum_{|\alpha| \leq 10} \int_0^{T_\varepsilon} \|Z^\alpha Q(u'_{k-1})(s, \cdot)\|_{L^2(\mathbb{R}^3)} \, ds.$$

Clearly for $|\alpha| \leq 10$ we have the pointwise bound

$$|Z^\alpha Q(u'_{k-1})(s,x)| \leq C_1 \left( \sum_{|\alpha| \leq 10} |Z^\alpha u'_{k-1}(s,x)| \right) \left( \sum_{|\alpha| \leq 5} |Z^\alpha u'_{k-1}(s,x)| \right).$$

Therefore, if we use (2.7), we have for a given $R = 2^j$, $j \geq 0$,

$$\sum_{|\alpha| \leq 10} \|Z^\alpha Q(u'_{k-1})(s, \cdot)\|_{L^2(\{|x| \in [2^j, 2^{j+1})\})}$$

$$\leq C 2^{-j} \sum_{|\alpha| \leq 10} \|Z^\alpha u'_{k-1}(s, \cdot)\|_{L^2(\{|x| \in [2^j, 2^{j+1})\})} \sum_{|\alpha| \leq 7} \|Z^\alpha u'_{k-1}(s, \cdot)\|_{L^2(\{|x| \in [2^{j-1}, 2^{j+2})\})}$$

$$\leq C \sum_{|\alpha| \leq 10} \|(1+r)^{-1/2} Z^\alpha u'_{k-1}(s, \cdot)\|^2_{L^2(\{|x| \in [2^{j-1}, 2^{j+2})\})}.$$

We also of course have the bounds

$$\sum_{|\alpha| \leq 10} \|Z^\alpha Q(u'_{k-1})(s, \cdot)\|_{L^2(\{|x| < 1\})} \leq C \sum_{|\alpha| \leq 10} \|Z^\alpha u'_{k-1}(s, \cdot)\|^2_{L^2(\{|x| < 2\})},$$

just by the standard Sobolev lemma. Thus, if we sum over $j$ and use (3.4) we get that

$$M_k(T_\varepsilon) \leq C_0 \varepsilon + C \sum_{|\alpha| \leq 10} \int_0^{T_\varepsilon} \|(1+r)^{-1/2} Z^\alpha u'_{k-1}(s, \cdot)\|^2_2 \, ds$$

$$= C_0 \varepsilon + C \sum_{|\alpha| \leq 10} \|(1+r)^{-1/2} Z^\alpha u'_{k-1}\|^2_{L^2(\{(s,x): 0 \leq s \leq T_\varepsilon\})}$$

$$\leq C_0 \varepsilon + C \ln(2 + T_\varepsilon) M^2_{k-1}(T_\varepsilon)$$

$$\leq C_0 \varepsilon + C \ln(2 + T_\varepsilon)(2 C_0 \varepsilon)^2,$$

using the induction hypothesis in the last step. If $\varepsilon$ and the $c$ occurring in the definition of $T_\varepsilon$ are small enough so that $4 C C_0 \varepsilon \ln(2 + T_\varepsilon) < 1$, then we get (3.3).

To show that the sequence $u_k$ converges to a solution, we estimate the quantity

$$(3.5) \quad A_k(T) = \sup_{0 \leq t \leq T} \sum_{|\alpha| \leq 10} \Big( \|Z^\alpha(u'_k - u'_{k-1})(t, \cdot)\|_{L^2(\mathbb{R}^3)}$$

$$+ (\ln(2 + t))^{-1/2} \|(1+r)^{-1/2} Z^\alpha(u'_k - u'_{k-1})\|_{L^2(\{(s,x): 0 \leq s \leq t\})} \Big).$$

The proof of Theorem 2.3 is concluded by the bounds

$$(3.6) \quad A_k(T_\varepsilon) \leq \frac{1}{2} A_{k-1}(T_\varepsilon), \quad k = 1, 2, \ldots.$$

Since $G$ is quadratic, we may repeat the above arguments to obtain

$$A_k(T_\varepsilon) \leq C \sum_{|\alpha| \leq 10} \int_0^{T_\varepsilon} \|(1+r)^{-1/2} Z^\alpha(u'_{k-1} - u'_{k-2})(s, \cdot)\|_{L^2(\mathbb{R}^3)}$$

$$\times \Big( \|(1+r)^{-1/2} Z^\alpha u'_{k-1}(s, \cdot)\|_{L^2(\mathbb{R}^3)} + \|(1+r)^{-1/2} Z^\alpha u'_{k-2}(s, \cdot)\|_{L^2(\mathbb{R}^3)} \Big) ds.$$

By applying the Schwarz inequality we conclude that

$$A_k(T_\varepsilon) \leq C \ln(2 + T_\varepsilon) \big( M_{k-1}(T_\varepsilon) + M_{k-2}(T_\varepsilon) \big) \cdot A_{k-1}(T_\varepsilon),$$

which by (3.3) leads to (3.6) if $\varepsilon$ and $T_\varepsilon$ are as above. $\square$

4. **Almost global existence outside of nontrapping obstacles.**

Let $\mathcal{K} \subset \mathbb{R}^3$ be a nontrapping obstacle (see [8]). We wish to adapt the argument for Minkowski space given in the previous section to show that there are almost global solutions to (1.1). To do this, we first need to see that analogs of the two estimates used in the proof of Theorem 1.2 carry over to the nontrapping obstacle setting.

In what follows, we shall assume, as we may, that $\mathcal{K} \subset \{|x| < 1/2\}$.

Clearly the analog of (2.7) works here. Namely, if $h(x) = 0$, $x \in \partial\mathcal{K}$ then

$$(4.1) \qquad \|h\|_{L^\infty(R/2 \leq |x| \leq R)} \leq CR^{-1} \sum_{|\alpha| \leq 2} \|Z^\alpha h\|_{L^2(R/4 \leq |x| \leq 2R)}.$$

For $R > 2$, this estimate follows directly from (2.6). For $R \leq 2$ it follows just from standard Sobolev estimates.

The other ingredient that is needed is an analog of Theorem 2.3. For simplicity, so that we avoid issues regarding compatibility conditions, we shall just prove estimates for the inhomogeneous wave equation with zero initial data, since that is the only thing that is needed for our application. So, we shall prove a variation of (2.5) for the solution of

$$(4.2) \qquad \begin{cases} \Box w = F, & (t,x) \in \mathbb{R}_+ \times \mathbb{R}^3\backslash\mathcal{K} \\ w(t,x) = 0, & x \in \partial\mathcal{K} \\ w(0,x) = 0, & \partial_t w(0,x) = 0. \end{cases}$$

**Theorem 4.1.** *If $w$ is as in (4.2) then for any $N = 0, 1, 2, \ldots$*

$$(4.3)$$
$$\sum_{|\alpha| \leq N} \left( \|Z^\alpha w'(t,\cdot)\|_{L^2(\mathbb{R}^3\backslash\mathcal{K})} + (\ln(2+t))^{-1/2} \|(1+r)^{-1/2} Z^\alpha w'\|_{L^2(\{(s,x)\in[0,t]\times\mathbb{R}^3\backslash\mathcal{K}\})} \right)$$
$$\leq C \sum_{|\alpha| \leq N} \int_0^t \|Z^\alpha F(s,\cdot)\|_{L^2(\mathbb{R}^3\backslash\mathcal{K})}\, ds + C \sup_{0 \leq s \leq t} \sum_{|\alpha| \leq N-1} \|Z^\alpha F(s,\cdot)\|_{L^2(\mathbb{R}^3\backslash\mathcal{K})}$$
$$+ C \sum_{|\alpha| \leq N-1} \|Z^\alpha F\|_{L^2(\{(s,x)\in[0,t]\times\mathbb{R}^3\backslash\mathcal{K}\})}.$$

Let us first see that the estimate holds when the norm in the left is taken over $|x| < 2$. Clearly the first term in the left is under control since

$$\sum_{|\alpha|\leq N} \|Z^\alpha w'(t,\cdot)\|_{L^2(\{x\in\mathbb{R}^3\backslash\mathcal{K}:\, |x|<2\})} \leq C_N \sum_{|\alpha|\leq N} \|\partial_{t,x}^\alpha w'(t,\cdot)\|_{L^2(\mathbb{R}^3\backslash\mathcal{K})},$$

and standard arguments imply that the right hand side here is dominated by

$$(4.4) \quad \sum_{|\alpha|\leq N} \|\partial_{t,x}^\alpha w'(t,\cdot)\|_{L^2(\mathbb{R}^3\backslash\mathcal{K})}$$
$$\leq C \sum_{|\alpha|\leq N} \int_0^t \|\partial_{t,x}^\alpha F(s,\cdot)\|_{L^2(\mathbb{R}^3\backslash\mathcal{K})}\, ds + C \sum_{|\alpha|\leq N-1} \|\partial_{t,x}^\alpha F(t,\cdot)\|_{L^2(\mathbb{R}^3\backslash\mathcal{K})}.$$

Indeed, if $N = 0$, (4.4) is just the standard energy identity. To prove that (4.4) holds for $N$, assuming that it is valid when $N$ is replaced by $N-1$, one notes that since $\partial_t w$ vanishes on the boundary one has

$$\sum_{|\alpha| \leq N-1} \|\partial_{t,x}^\alpha \partial_t w'(t, \cdot)\|_{L^2(\mathbb{R}^3 \setminus \mathcal{K})}$$
$$\leq C \sum_{|\alpha| \leq N-1} \int_0^s \|\partial_{t,x}^\alpha \partial_s F(s, \cdot)\|_{L^2(\mathbb{R}^3 \setminus \mathcal{K})} + C \sum_{|\alpha| \leq N-2} \|\partial_{t,x}^\alpha \partial_t F(t, \cdot)\|_{L^2(\mathbb{R}^3 \setminus \mathcal{K})}.$$

Since $\partial_t^2 w = \Delta w + F$, we get from this that

$$\sum_{|\alpha| \leq N-1} \|\partial_{t,x}^\alpha \Delta w(t, \cdot)\|_{L^2(\mathbb{R}^3 \setminus \mathcal{K})}$$
$$\leq C \sum_{|\alpha| \leq N} \int_0^s \|\partial_{t,x}^\alpha F(s, \cdot)\|_{L^2(\mathbb{R}^3 \setminus \mathcal{K})} \, ds + \sum_{|\alpha| \leq N-1} \|\partial_x^\alpha F(t, \cdot)\|_{L^2(\mathbb{R}^3 \setminus \mathcal{K})}.$$

By elliptic regularity, $\sum_{|\alpha| \leq N} \|\partial_{t,x}^\alpha w'(t, \cdot)\|_{L^2(\mathbb{R}^3 \setminus \mathcal{K})}$ is dominated by the left side of the last equation, which finishes the proof of (4.4), since

$$\sum_{|\alpha| \leq N} \|\partial_{t,x}^\alpha w'(t, \cdot)\|_{L^2(\mathbb{R}^3 \setminus \mathcal{K})} \leq \sum_{|\alpha| \leq N} \|\partial_x^\alpha w'(t, \cdot)\|_{L^2(\mathbb{R}^3 \setminus \mathcal{K})} + \sum_{|\alpha| \leq N-1} \|\partial_{t,x}^\alpha \partial_t w'(t, \cdot)\|_{L^2(\mathbb{R}^3 \setminus \mathcal{K})}.$$

To handle the second term on the left side of (4.3), again when the left hand norm is taken over $|x| < 2$, we shall need the following

**Lemma 4.2.** *If $w$ is as in* (4.2) *then for any $N = 0, 1, 2, \ldots$*

$$(4.5) \quad \sum_{|\alpha| \leq N} \|\partial_{t,x}^\alpha w'\|_{L^2(\{(s,x) \in [0,t] \times \mathbb{R}^3 \setminus \mathcal{K}: |x| < 2\})}$$
$$\leq C \sum_{|\alpha| \leq N} \int_0^t \|\partial_{t,x}^\alpha F(s, \cdot)\|_{L^2(\mathbb{R}^3 \setminus \mathcal{K})} \, ds + C \sum_{|\alpha| \leq N-1} \|\partial_{t,x}^\alpha F\|_{L^2(\{(s,x) \in [0,t] \times \mathbb{R}^3 \setminus \mathcal{K}\})}.$$

Clearly (4.5) implies that

$$\sum_{|\alpha| \leq N} \|Z^\alpha w'\|_{L^2(\{(s,x) \in [0,t] \times \mathbb{R}^3 \setminus \mathcal{K}: |x| < 2\})}$$
$$\leq C \sum_{|\alpha| \leq N} \int_0^t \|Z^\alpha F(s, \cdot)\|_{L^2(\mathbb{R}^3 \setminus \mathcal{K})} \, ds + C \sum_{|\alpha| \leq N-1} \|Z^\alpha F\|_{L^2(\{(s,x) \in [0,t] \times \mathbb{R}^3 \setminus \mathcal{K}\})},$$

finishing the proof that the analog of (4.3) holds where the norms in the left are taken over $|x| < 2$.

**Proof of Lemma 4.2:** By the proof of (4.4), (4.5) follows from the special case where $N = 0$:

$$(4.6) \quad \|w'\|_{L^2(\{(s,x) \in [0,t] \times \mathbb{R}^3 \setminus \mathcal{K}: |x| < 2\})} \leq C \int_0^t \|F(s, \cdot)\|_{L^2(\mathbb{R}^3 \setminus \mathcal{K})} \, ds.$$

Let us first prove (4.6) under the assumption that $F(s,x) = 0$ when $|x| > 4$. In this case the estimate just follows from the fact that there is a $c > 0$ so that

$$\|w'(t,\,\cdot\,)\|_{L^2(\{x\in\mathbb{R}^3\setminus\mathcal{K}:\,|x|<4\})} \leq C\int_0^t e^{-c(t-s)}\|F(s\,\cdot\,)\|_2\,ds, \quad \text{if } F(s,x) = 0,\ |x| > 4,$$

because of the exponential local energy decay of Morawetz, Ralston, and Strauss [8] (see also [6]), and our assumption that the obstacle $\mathcal{K}$ is nontrapping.

Note that this estimate also implies that if $v$ solves the Dirichlet-wave equation, $\Box v = G$, $v(t,x) = 0$, $x \in \partial\mathcal{K}$ with zero initial data then

$$\|v'\|_{L^2(\{(s,x)\in[0,t]\times\mathbb{R}^3\setminus\mathcal{K}:\,|x|<4\})} \leq C\|G\|_{L^2(\{(s,x)\in[0,t]\times\mathbb{R}^3\setminus\mathcal{K}:\,|x|<4\})},$$
$$\text{if } G(s,x) = 0,\ |x| > 4.$$

Using this we can prove (4.6) for the remaining case where we assume that $F(s,x) = 0$, $|x| \leq 4$. To do this we choose $\rho \in \mathbb{R}^3$ satisfying $\rho(x) = 1$, $|x| \leq 2$, and $\rho(x) = 0$, $|x| > 4$. If we then write $w = u_0 + u_r$ where $u_0$ solves the boundaryless inhomogeneous wave equation $\Box u_0 = F$ with zero initial data, we then set $v = \rho u_0 + u_r$ and note that $\Box v = -2\nabla_x\rho \cdot \nabla_x u_0 - (\Delta\rho)u_0$ vanishes for $|x| < 4$. Therefore,

$$\|w'\|_{L^2(\{(s,x)\in[0,t]\times\mathbb{R}^3\setminus\mathcal{K}:\,|x|<2\})} = \|v'\|_{L^2(\{(s,x)\in[0,t]\times\mathbb{R}^3\setminus\mathcal{K}:\,|x|<2\})}$$
$$\leq C\|u_0'\|_{L^2(\{(s,x)\in[0,t]\times\mathbb{R}^3\setminus\mathcal{K}:\,|x|<4\})} + C\|u_0\|_{L^2(\{(s,x)\in[0,t]\times\mathbb{R}^3\setminus\mathcal{K}:\,|x|<4\})}.$$

By (2.2) and (2.5) the right side here is dominated by that of (4.6), which finishes the proof. $\square$

**End of proof of Theorem 4.1:** We need to see that

(4.7)
$$\sum_{|\alpha|\leq N}\left(\|Z^\alpha w'(t,\,\cdot\,)\|_{L^2(|x|>2)} + (\ln(2+t))^{-1/2}\|(1+r)^{-1/2}Z^\alpha w'\|_{L^2(\{[0,t]\times\{x:\,|x|>2\}\})}\right)$$
$$\leq C\sum_{|\alpha|\leq N}\int_0^t \|Z^\alpha F(s,\,\cdot\,)\|_{L^2(\mathbb{R}^3\setminus\mathcal{K})}\,ds + C\sup_{0\leq s\leq t}\sum_{|\alpha|\leq N-1}\|Z^\alpha F(s,\,\cdot\,)\|_{L^2(\mathbb{R}^3\setminus\mathcal{K})}$$
$$+ C\sum_{|\alpha|\leq N-1}\|Z^\alpha F\|_{L^2(\{(s,x)\in[0,t]\times\mathbb{R}^3\setminus\mathcal{K}\})}.$$

For this we fix $\beta \in C^\infty(\mathbb{R}^3)$ satisfying $\beta(x) = 1$, $|x| \geq 2$ and $\beta(x) = 0$, $|x| \leq 1$. Then since, by assumption the obstacle is contained in the set $|x| < 1/2$, it follows that $v = \beta w$ solves the boundaryless wave equation

$$\Box v = \beta F - 2\nabla_x\beta \cdot \nabla_x w - (\Delta\beta)w$$

with zero initial data, and satisfies $w(t,x) = v(t,x)$, $|x| \geq 2$. If we split $v = v_1 + v_2$, where $\Box v_1 = \beta F$, and $\Box v_2 = -2\nabla_x\beta \cdot \nabla_x w - (\Delta\beta)w$, it then suffices to prove that

(4.8)
$$\sum_{|\alpha|\leq N} \left( \|Z^\alpha v_2'(t,\,\cdot\,)\|_{L^2(|x|>2)} + (\ln(2+t))^{-1/2}\|(1+r)^{-1/2}Z^\alpha v_2'\|_{L^2(\{[0,t]\times\{x:\,|x|>2\}\})} \right)$$
$$\leq C \sum_{|\alpha|\leq N} \int_0^t \|Z^\alpha F(s,\,\cdot\,)\|_{L^2(\mathbb{R}^3\setminus\mathcal{K})}\,ds + C \sup_{0\leq s\leq t} \sum_{|\alpha|\leq N-1} \|Z^\alpha F(s,\,\cdot\,)\|_{L^2(\mathbb{R}^3\setminus\mathcal{K})}$$
$$+ C \sum_{|\alpha|\leq N-1} \|Z^\alpha F\|_{L^2(\{(s,x)\in[0,t]\times\mathbb{R}^3\setminus\mathcal{K}\})}.$$

This is because by (2.6) we have

$$\sum_{|\alpha|\leq N} \left( \|Z^\alpha v_1'(t,\,\cdot\,)\|_{L^2(|x|>2)} + (\ln(2+t))^{-1/2}\|(1+r)^{-1/2}Z^\alpha v_1'\|_{L^2(\{[0,t]\times\{x:\,|x|>2\}\})} \right)$$
$$\leq C \sum_{|\alpha|\leq N} \int_0^t \|Z^\alpha F(s,\,\cdot\,)\|_2\,ds,$$

due to the fact that

$$\sum_{|\alpha|\leq N} \int_0^t \|Z^\alpha(\beta F)(s,\,\cdot\,)\|_2\,ds \leq C \sum_{|\alpha|\leq N} \int_0^t \|Z^\alpha F(s,\,\cdot\,)\|_2\,ds.$$

To prove (4.8) we note that $G = -2\nabla_x\beta \cdot \nabla_x w - (\Delta\beta)w = \Box v_2$, vanishes unless $1 < |x| < 2$. To use this, fix $\chi \in C_0^\infty(\mathbb{R})$ satisfying $\chi(s) = 0$, $|s| > 2$, and $\sum_j \chi(s-j) = 1$. We then split $G = \sum_j G_j$, where $G_j(s,x) = \chi(s-j)G(s,x)$, and let $v_{2,j}$ be the solution of the corresponding inhomogeneous wave equation $\Box v_{2,j} = G_j$ with zero initial data in Minkowski space. By sharp Huygen's principle we have that $|Z^\alpha v_2(t,x)|^2 \leq C \sum_j |Z^\alpha v_{2,j}(t,x)|^2$ for some uniform constant $C$. Therefore, by (2.6) we have that the square of the left side of (4.8) is dominated by

$$\sum_{|\alpha|\leq N} \sum_j \left( \int_0^t \|Z^\alpha G_j(s,\,\cdot\,)\|_2\,ds \right)^2$$
$$\leq C \sum_{|\alpha|\leq N} \|Z^\alpha G\|^2_{L^2(\{(s,x):\,0\leq s\leq t,\,1<|x|<2\})}$$
$$\leq C \sum_{|\alpha|\leq N} \|Z^\alpha w'\|^2_{L^2(\{(s,x):\,0\leq s\leq t,\,1<|x|<2\})} + C \sum_{|\alpha|\leq N} \|Z^\alpha w\|^2_{L^2(\{(s,x):\,0\leq s\leq t,\,1<|x|<2\})}$$
$$\leq C \sum_{|\alpha|\leq N} \|Z^\alpha w'\|^2_{L^2(\{(s,x)\in[0,t]\times\mathbb{R}^3\setminus\mathcal{K}:\,|x|<2\})}$$
$$\leq C \sum_{|\alpha|\leq N} \|\partial_{t,x}^\alpha w'\|^2_{L^2(\{(s,x)\in[0,t]\times\mathbb{R}^3\setminus\mathcal{K}:\,|x|<2\})}.$$

Consequently, (4.8) follows from (4.5), which finishes the proof. $\square$

Let us conclude by showing how we can adapt the proof of Theorem 1.2 to show that there is almost global existence for (1.1). As in [3] it is convenient to show that one can solve an equivalent nonlinear equation which has zero initial data to avoid having to deal with issues regarding compatibility conditions for the data. We can then set up an iteration argument for this new equation that is similar to the one used in the proof of Theorem 1.2.

To make the reduction we first note that if the data satisfies (1.3) with $\varepsilon$ small we can find a local solution $u$ to $\Box u = Q(u')$ in $0 < t < 1$. We then claim that if (1.3) is valid, with $\varepsilon > 0$ sufficiently small then there must be a constant $C$ so that

$$(4.9) \quad \sup_{0 \leq t \leq 1} \sum_{|\alpha| \leq 10} \Big( \|Z^\alpha u'(t, \,\cdot\,)\|_{L^2(\mathbb{R}^3 \setminus \mathcal{K})} + \|(1+r)^{-1/2} Z^\alpha u'\|_{L^2(\{(s,x) \in [0,1] \times \mathbb{R}^3 \setminus \mathcal{K}\})} \Big) \leq C\varepsilon.$$

One sees this by first noticing that standard local existence theory (see e.g., [3]) implies that if $\varepsilon$ is small then the analog of (4.9) holds if the norms on the left side are taken over, say, the region where $|x| < 10$. But then (3.1) implies that the estimate where the norms are taken over $|x| > 10$ is also valid by the finite propagation speed for $\Box$ since in this region if $0 < t < 1$, $u$ agrees with a solution of the boundaryless nonlinear wave equation $\Box u = Q(u')$ with data equal to a cutoff function times the original data $(f, g)$ for (1.1).

Using this local solution we can set up our iteration. We first fix a bump function $\eta \in C^\infty(\mathbb{R})$ satisfying $\eta(t) = 1$ if $t \leq 1/2$ and $\eta(t) = 0$ if $t > 1$. If we set

$$u_0 = \eta u$$

then

$$\Box u_0 = \eta Q(u') + [\Box, \eta] u.$$

So $u$ will solve $\Box u = Q(u')$ for $0 < t < T_\varepsilon$ if and only if $w = u - u_0$ solves

$$(4.10) \quad \begin{cases} \Box w = (1-\eta) Q((u_0+w)') - [\Box, \eta](u_0+w) \\ w|_{\partial \mathcal{K}} = 0 \\ w(0, x) = \partial_t w(0, x) = 0 \end{cases}$$

for $0 < t < T_\varepsilon$.

We shall solve this equation by iteration. We set $w_0 = 0$ and then define $w_k$, $k = 1, 2, 3, \ldots$ inductively by requiring that

$$(4.11) \quad \begin{cases} \Box w_k = (1-\eta) Q((u_0 + w_{k-1})') - [\Box, \eta](u_0 + w_k) \\ w_k|_{\partial B} = 0 \\ w_k(0, x) = \partial_t w_k(0, x) = 0. \end{cases}$$

As before, we let

$$M_k(T) = \sup_{0 \leq t \leq T} \sum_{|\alpha| \leq 10} \Big( \|Z^\alpha w_k'(t, \,\cdot\,)\|_2 + (\ln(2+t))^{-1/2} \|(1+r)^{-1/2} Z^\alpha w_k'\|_{L^2(\{(s,x):\, 0 \leq s \leq t\})} \Big).$$

Then since $u_0$ vanishes for $t > 1$ and satisfies the bounds in (4.9), we can argue as in §3, using now (4.1) and (4.3) to see that if $\varepsilon < \varepsilon_0$ is sufficiently small and if the $c$ occurring in the definition of $T_\varepsilon$ is sufficiently small then for every $k = 1, 2, 3, \ldots$ we have

$$(4.12) \qquad M_k(T_\varepsilon) \leq C_0 \varepsilon,$$

for some uniform constant $C_0$. The additional two terms in the right side of (4.3) as opposed to (2.6) do not cause any problems for $\varepsilon > 0$ small. Indeed, since every occurrence of $u_0$ contributes $O(\varepsilon)$ by (4.9), the second term in the right side of (4.3) will contribute $\leq C(\varepsilon + M_{k-1}(T_\varepsilon))^2$ to the bounds for $M_k(T_\varepsilon)$, while the last term in the right will contribute $\leq C \log(2 + T_\varepsilon)(\varepsilon + M_{k-1}(T_\varepsilon))^2$, which means that the arguments from §3 will lead to the bounds

$$M_k(T_\varepsilon) \leq C_1 \varepsilon + C_1 \ln(2 + T_\varepsilon)(\varepsilon + M_{k-1}(T_\varepsilon))^2 + C_1(\varepsilon + M_{k-1}(T_\varepsilon))^2,$$

for some uniform constant $C_1$, if $\varepsilon$ is small, which implies by induction that (4.12) holds with $C_0 = 2C_1$ for such $\varepsilon$ since $M_0 \equiv 0$.

Furthermore, if we let $A_k$ be the analog of (3.5) then the same arguments as in §3 will show that (3.6) holds for small enough $\varepsilon$. Because of this and (4.12) we conclude that $w_k$ must converge to a solution of (4.10) which satisfies the analog of (3.1). This means that $u = u_0 + w$ will be a solution of our original equation (1.1) also verifying the analog of (3.1). If the data is $C^\infty$ and satisfies the compatibility conditions to infinite order the solution will be $C^\infty$ on $[0, T_\varepsilon] \times \mathbb{R}^3 \setminus \mathcal{K}$ by standard local existence theory (see e.g., [3]). This completes the proof of Theorem 1.1.

Department of Mathematics, University of Minnesota, Minnesota, MN 55455

Department of Mathematics, University of Washington, Seattle, WA 98195

Department of Mathematics, The Johns Hopkins University, Baltimore, MD 21218